\documentclass{amsart}

\usepackage[all]{xy}
\usepackage{amscd}
\usepackage[mathscr]{eucal}
\usepackage{amssymb}

\newcommand{\colim}{\operatornamewithlimits{colim}}
\DeclareMathOperator{\Map}{Map}
\newtheorem{lem}{Lemma}[section]
\newtheorem{defi}[lem]{Definition}
\newtheorem{prop}[lem]{Proposition}
\newtheorem{rema}[lem]{Remark}
\newtheorem{theo}[lem]{Theorem}
\newtheorem{coro}[lem]{Corollary}
\newtheorem*{exam}{Example}
\newtheorem*{theorem}{Theorem}
\newtheorem*{corollary}{Corollary}

\begin{document}
\title{Splitting With Continuous Control in Algebraic $K$-theory}
\author{David Rosenthal}
\address{Department of Mathematics and Statistics, McMaster University, Hamilton, ON L8S 4K1, Canada}
\email{rosend@math.mcmaster.ca}

\begin{abstract}
	In this work, the continuously controlled assembly map in algebraic $K$-theory, as developed by Carlsson and Pedersen, is proved to be a split injection for groups $\Gamma$ that satisfy certain geometric conditions. The group $\Gamma$ is allowed to have torsion, generalizing a result of Carlsson and Pedersen. Combining this with a result of John Moody, $K_0(k\Gamma)$ is proved to be isomorphic to the colimit of $K_0(kH)$ over the finite subgroups $H$ of $\Gamma$, when $\Gamma$ is a virtually polycyclic group and $k$ is a field of characteristic zero.
\end{abstract}

\maketitle

\section{Introduction}
	The algebraic $K$-theory groups $K_n(R\Gamma)$ contain geometric information about manifolds with fundamental group $\Gamma$. One way to explore these groups is to study an assembly map that relates the algebraic $K$-theory of the group ring $R\Gamma$ to the homology of $\Gamma$. This is formulated as a map from a generalized homology theory evaluated on the universal space for a chosen family of subgroups of $\Gamma$, to the algebraic $K$-theory of $R\Gamma$. The goal here is to show that for the family of finite subgroups, the assembly map is a {\it split injection} for groups whose universal space satisfies certain geometric conditions. There are many ways to build assembly maps. The following theorem is proved in this paper using the {\it continuously controlled} model developed by Carlsson and Pedersen~\cite{cp}.

\begin{theorem}[Theorem~\ref{main}]
	Let $\Gamma$ be a discrete group, and let $\mathcal{E}=E\Gamma(\mathfrak{f})$ be the universal space for $\Gamma$-actions with finite isotropy, where $\mathfrak{f}$ denotes the family of finite subgroups of $\Gamma$. Assume that $\mathcal{E}$ is a finite $\Gamma$-CW complex admitting a compactification, $X$, (i.e., $X$ is compact and $\mathcal{E}$ is an open dense subset) such that
	\begin{enumerate}
		\item[1.] the $\Gamma$-action extends to $X$;
		\item[2.] $X$ is metrizable;
		\item[3.] $X^G$ is contractible for every $G\in \mathfrak{f}$;
		\item[4.] ${\mathcal{E}}^G$ is dense in $X^G$ for every $G\in \mathfrak{f}$;
		\item[5.] compact subsets of $\mathcal{E}$ become small near $Y=X-\mathcal{E}$.  That is, for every compact subset $K \subset \mathcal{E}$ and for every neighborhood $U \subset X$ of $y \in Y$, there exists a neighborhood $V \subset X$ of $y$ such that $g \in \Gamma$ and $gK \cap V \neq \emptyset$ implies $gK \subset U$.
	\end{enumerate}
Then $H_*^{\Gamma}(\mathcal{E};{\mathbb{K}}^{-\infty}(R{\Gamma}_x)) \to K_*(R\Gamma)$ is split injective for every ring $R$.
\end{theorem}

If $\Gamma$ is torsion-free, then $E\Gamma(\mathfrak{f})=E\Gamma$, and Theorem~\ref{main} is~\cite[Theorem A]{cp}. In~\cite{dl}, Davis-L\"{u}ck provide a different formulation of the assembly map. However, in recent work of Hambleton-Pedersen~\cite{hp}, it is shown that the Davis-L\"{u}ck version is homotopy equivalent to the continuously controlled version. The continuously controlled model of the assembly map that is used in~\cite{hp} is actually different from the one presented here. However, they too are homotopy equivalent. The benefit of using the version from \cite{cp} is that it realizes the assembly map as a map of fixed sets. This admits the possibility of proving that it induces a split injection on homotopy groups with the use of homotopy fixed sets. Specifically we construct an equivariant map of spectra $S \to T$, as in \cite{cp}, such that the induced map, $\pi_*(S^{\Gamma}) \to \pi_*(T^{\Gamma})$, is $H_*^{\Gamma}(\mathcal{E};{\mathbb{K}}^{-\infty}(R{\Gamma}_x)) \to K_*(R\Gamma)$. Then the diagram
\[ \xymatrix{
	S^{\Gamma} \ar[d]_{g} \ar[r]^{f} & T^{\Gamma} \ar[d] \\
	S^{h\mathfrak{f} \Gamma} \ar[r]^{h} & T^{h\mathfrak{f} \Gamma} } \]
is examined, where $S^{h\mathfrak{f} \Gamma}$ denotes the set of all $\Gamma$-equivariant maps from $E\Gamma(\mathfrak{f})$ to $S$. The diagram is induced by the $\Gamma$-equivariant map, $E\Gamma(\mathfrak{f}) \to \bullet$. If $\pi_*(g)$ and $\pi_*(h)$ are shown to be isomorphisms, then $\pi_*(f)$ is a split injection. That is, there is a homomorphism $\phi:\pi_*(T^{\Gamma}) \to \pi_*(S^{\Gamma})$ such that $\phi \circ \pi_*(f)=1_{\pi_*(S^{\Gamma} )}$. This would imply that $\pi_*(S^{\Gamma})$ is a direct summand of $\pi_*(T^{\Gamma})\cong K_*(R\Gamma)$. 

	A result of Wilking~\cite[Theorem 3]{wilking} shows that virtually polycyclic groups satisfy the conditions of Theorem~\ref{main}, since their universal space can be chosen to be some Euclidean $n$-space. The desired compactification is then obtained by compactifying the space with an $(n-1)$-sphere at infinity. Therefore the assembly map described above is a split injection for virtually polycyclic groups. In this case the assembly map in dimension zero can be identified with the induction map, $\colim_{H \in \mathfrak{f}} K_0(kH) \to K_0(k\Gamma)$, where $k$ is a field of characteristic zero (Proposition~\ref{h_0}). A theorem of Moody~\cite[Theorem 1]{moody} shows that this map is a surjection. This yields the following corollary.

	\begin{corollary}[Corollary~\ref{iso}]
	Let $\Gamma$ be a virtually polycyclic group and $k$ a field of characteristic zero. Then the induction map
\[ \colim_{H \in \mathfrak{f}} K_0(kH) \to K_0(k\Gamma) \]
is an isomorphism.
	\end{corollary}

	Sections~\ref{cca} and~\ref{hcca} contain the material from continuously controlled algebra required for the proof of Theorem~\ref{main}. In particular, Section~\ref{hcca} connects continuosly controlled algebra to two functors that give rise to two generalized homology theories, one of which is a generalized reduced Steenrod homology theory. For a more thorough treatment of the theory of continuously controlled algebra the reader is encouraged to read~\cite{cp}. The homotopy theory necessary to use the generalized fixed sets and the diagram described above is provided in Section~\ref{gfs}. Lemma~\ref{lemmaa} in this section makes it clear that in order to prove that the map $h$ in the above diagram is a weak homotopy equivalence, it will be necessary to work with finite groups. It turns out that it will also be necessary to show that the reduced Steenrod homology (built in Section~\ref{hcca}) of the quotient of a contractible compact metrizable space by a finite group is trivial. This is proved in Section~\ref{ccvsh}, making use of the Conner Conjecture.

	In Sections~\ref{akcca} and~\ref{fccfs}, Theorem~\ref{main} is proved. Section~\ref{akcca} contains a proof that the map $g$ in the above diagram is a weak homotopy equivalence. The proof is adapted from~\cite{cp}. In Section~\ref{fccfs}, the most difficult step in the proof of Theorem~\ref{main} is handled. Here the map $h$ in the above diagram is proved to be a weak homotopy equivalence using a filtration of the space $X$ by {\it conjugacy classes of fixed sets}. In Section~\ref{dlss} a spectral sequence from the work of Davis-L\"{u}ck~\cite{dl} is developed. This is then used in Section~\ref{vpg} to show that Theorem~\ref{main} can be combined with Moody's result to prove Corollary~\ref{iso}. Some of the results in the paper are taken from the author's doctoral dissertation~\cite{me}.

	I would like to thank my thesis supervisor, Erik K. Pedersen, for all of his advice and positive energy. I want to thank Jim Davis for suggesting Corollary~\ref{iso} and for pointing me in the right direction. I would also like to thank Ian Hambleton and Andrew Nicas for many useful conversations during the preparation of this paper, in particular with the material in Sections~\ref{dlss} and~\ref{vpg}.

\section{Continuously Controlled Algebra}\label{cca}

	Let $\Gamma$ be a discrete group acting on a topological space $E$ with finite isotropy, and let $R$ be a ring with unit. The free $R$-module generated by $E \times \Gamma \times \mathbb{N}$ is denoted $R[E \times \Gamma]^{\infty}$. If $A$ is a submodule of $R[E \times \Gamma]^{\infty}$, denote $A \cap R[x \times \Gamma]^{\infty}$ by $A_x$, where $x \in E$.

	Define the small additive category $\mathcal{B}(E;R)$ to have as objects submodules of $R[E \times \Gamma]^{\infty}$, A, where
\begin{enumerate}
	\item[(i)] $A=\bigoplus_{x \in E} A_x$;
	\item[(ii)] $A_x$ is a finitely generated free $R$-module with basis contained in $\{ x \} \times \mathbb{N}$;
	\item[(iii)] $\{x \in E \, | \, A_x \neq 0 \}$ is locally finite in $E$.
\end{enumerate}
Morphisms are all $R$-module morphisms. Given a morphism $\phi:A\to B$, and two points $x,y\in E$, we denote the {\it component} of $\phi$ from $A_x$ to $B_y$ to be $\phi^x_y$. Thus, $\phi$ is determined by, and determines a collection of $R$-homomorphisms, $\{\phi^x_y:A_x \to B_y \}$.

	The diagonal action of $\Gamma$ on $E \times \Gamma$ provides a $\Gamma$-action on $R[E \times \Gamma]^{\infty}$. This induces a $\Gamma$-action on $\mathcal{B}(E;R)$ such that the fixed category $\mathcal{B}(E;R)^{\Gamma}$, has those objects, $A$, in $\mathcal{B}(E;R)$, that satisfy the conditions
\begin{enumerate}
  \item[1.] $A_{\gamma x}\cong A_x$ for every $\gamma \in \Gamma$, and
  \item[2.] $A_x$ is a finitely generated free $R{\Gamma}_x$-module. 
\end{enumerate}  
This implies that $\bigoplus_{x'\in [x]} A_{x'}$ is a finitely generated free $R\Gamma$-module, where $[x]=\{ \gamma x \, | \, \gamma \in \Gamma \}$. The morphisms in $\mathcal{B}(E;R)^{\Gamma}$ are those morphisms, $\phi$, in $\mathcal{B}(E;R)$ that are $\Gamma$-equivariant, i.e., $\gamma {\phi}^x_y {\gamma}^{-1}=\phi^{\gamma x}_{\gamma y}$ for all $\gamma \in \Gamma$ and all $x,y \in E$.

	This definition of $\mathcal{B}(E;R)$ differs from the one given in~\cite{cp}. Here objects come equipped with a basis, but this does not affect the $K$-theory. Objects are required to have a basis in $E \times \Gamma \times \mathbb{N}$ so that those that are in the fixed category have the properties listed above which are needed to handle groups with torsion. Crossing with $\Gamma$ is also required in order to do this. In~\cite{cp}, objects are submodules of $R[E]^{\infty}$ instead of $R[E \times \Gamma]^{\infty}$. The $\Gamma$-actions on these two $R$-modules induce different actions on $\mathcal{B}(E;R)$, thus producing slightly different fixed categories. However, if $\Gamma$ is a torsion free group, as in~\cite{cp}, the two fixed categories agree (modulo the basis requirement). Furthermore, it is easy to see that the results from~\cite[Section 1]{cp} hold with the definition used here.
	
	Let us recall some categorical terminology.  Two functors, $F$ and $G$, between categories $\mathcal{A}$ and $\mathcal{B}$, are {\it naturally equivalent} if there is a natural transformation from $F$ to $G$ that is an isomorphism for every object in $\mathcal{A}$.  Categories $\mathcal{A}$ and $\mathcal{B}$ are {\it equivalent} if there are functors, $F:\mathcal{A} \to \mathcal{B}$ and $G:\mathcal{B} \to \mathcal{A}$, such that $FG$ is naturally equivalent to $1_{\mathcal{B}}$ and $GF$ is naturally equivalent to $1_{\mathcal{A}}$. Two additive categories are {\it isomorphic} when they are equivalent by functors that give a one-one correspondence between objects. The symbol $=$ is used to denote an isomorphism 
that is the identity on objects.

	Let $X$ be a topological space, $Y$ a closed subspace of $X$, and $E=X-Y$.

\begin{defi}
{\rm A morphism, $\phi:A \to B$ in $\mathcal{B}(E;R)$ is {\it continuously controlled at} $y\in Y$ if for every neighborhood $U \subseteq X$ of $y$, there exists a neighborhood $V \subseteq X$ of $y$ such that $\phi^a_b =0$ and $\phi^b_a =0$ whenever $a \in V$ and $b \notin U$.}
\end{defi}

The category $\mathcal{B}(X,Y;R)$ has the same objects as $\mathcal{B}(E;R)$, but the morphisms are required to be continuously controlled at every point in $Y$.

	We require several categories related to $\mathcal{B}(X,Y;R)$. Let $W$ be an open subset of $Y$, $C$ a closed subset of $Y$, $T$ an open subset of $X$, $K$ a topological space, and $p:T\to K$ a function that is continuous on $Y\cap T$.

\begin{defi}
{\rm A morphism $\phi:A\to B$ in $\mathcal{B}(E;R)$ is said to be {\it $p$-controlled at} $y\in Y \cap T$, if for every neighborhood $U \subseteq K$ of $p(y)$, there is a neighborhood $V \subseteq X$ of $y$ such that $\phi^a_b =0$ and $\phi^b_a =0$ whenever $a \in V$ and $b \notin p^{-1}(U)$.}
\end{defi}

The category $\mathcal{B}(X,Y,p;R)$ has the same objects as $\mathcal{B}(E;R)$, but the morphisms are required to be continuously controlled at all points of $Y-T$ and $p$-controlled at all points of $Y \cap T$.
There are only two categories of this type that we will use: $\mathcal{B} (CX,CY \cup X,p_X;R)$ and $\mathcal{B} (\Sigma X,\Sigma Y,p_X;R)$, where $CX$ denotes the cone of $X$, $\Sigma X$ denotes the unreduced suspension of $X$, and $p_X: X\times (0,1) \to X$ is the projection map.

The category $\mathcal{B}(X,Y;R)^W$ has the same objects as $\mathcal{B}(X,Y;R)$, but morphisms are identified if they agree in a neighborhood of $W$. Specifically, $\phi, \psi:A\to B$ are identified if there is a neighborhood $U \subseteq X$ of $W$ so that $a\in U$ or $b\in U$ implies that ${\phi}^a_b = {\psi}^a_b$. Similarly, the category $\mathcal{B}(X,Y,p;R)^W$ has the same objects as $\mathcal{B}(X,Y,p;R)$, but morphisms are identified if they agree on a neighborhood of $W$.  We will refer to $\mathcal{B}(X,Y;R)^W$ as a {\it germ} category.

\begin{defi}
{\rm The {\it support at infinity} of an object $A$ in $\mathcal{B}(X,Y;R)$, denoted $supp_{\infty}(A)$, is the set of limit points of $\{x \, | \, A_x \neq 0\}$.}
\end{defi}

The category $\mathcal{B}(X,Y;R)_C$ is the full subcategory of $\mathcal{B}(X,Y;R)$ on objects, $A$, with $supp_{\infty}(A)\subseteq C$. Similarly, $\mathcal{B}(X,Y;R)_C^W$, $\mathcal{B}(X,Y,p;R)_C$, and $\mathcal{B}(X,Y,p;R)_C^W$ are full subcategories of $\mathcal{B}(X,Y;R)^W$, $\mathcal{B}(X,Y,p;R)$, and $\mathcal{B}(X,Y,p;R)^W$ respectively.

\begin{defi}
{\rm	A function $f:(X,Y) \to (X',Y')$ is {\it eventually continuous} if
\begin{enumerate}
	\item[(i)] for every compact $K \subset X'-Y'$, $f^{-1}(K)$ has compact closure in $X-Y$;
	\item[(ii)] $f(X-Y) \subset X'-Y'$;
	\item[(iii)] $f$ is continuous on $Y$.
\end{enumerate}}
\end{defi}

\begin{prop}\cite[Proposition 1.16]{cp}\label{eventually}
	An eventually continuous map $f:(X,Y) \to (X',Y')$ induces a functor $f_{\sharp}:\mathcal{B}(X,Y;R) \to \mathcal{B}(X',Y';R)$, defined as follows. An object $A=\{A_x\}$ is sent to $f_{\sharp}A$, where $(f_{\sharp}A)_{x'}=\bigoplus^{}_{x\in f^{-1}(x')} A_x$, and morphisms are induced by the identity.
If $f$ sends $Y-W$ to $Y'-W'$, then $f$ also induces a functor $f_{\sharp}:\mathcal{B}(X,Y;R)^W \to \mathcal{B}(X',Y';R)^{W'}$.
\end{prop}

\begin{prop}\cite[Proposition 1.18]{cp}\label{equal}
	If $f_1$, $f_2:(X,Y) \to (X',Y')$ are eventually continuous maps and $f_1|_Y=f_2|_Y$, then $f_{1\sharp}$ and $f_{2\sharp}$ are naturally equivalent functors.
\end{prop}

        The next two lemmas build up to Lemma~\ref{p_x} which will be used in Section~\ref{fccfs}.

\begin{lem}\label{retract}
	Let $X'$ be a compact metrizable space, let $Y'$ be a closed nowhere dense subset of $X'$, and let $Z'$ be a closed subset of $X'$ such that $\overline{Z'-Y'}=Z'$. Then $\mathcal{B}(X',Y';R)_{Z' \cap Y'} \cong \mathcal{B}(Z',Z' \cap Y';R)$.
\end{lem}

\begin{proof}
	Let $x\in X'-(Z'\cup Y')$ be given. Choose $y_x\in Z'\cap Y'$ so that $d(x,y_x)=d(x,Z'\cap Y')$. Since $Z'-Y'$ is dense in $Z'$, we can choose $z_x\in Z'-Y'$ so that $d(z_x,y_x)\leq d(x,Z'\cap Y')$. Define $f:(X',Y')\to (X',Y')$ by
\[ f(x)=\left\{
\begin{array}{cl}
	z_x & {\rm if} \ x \in X'-(Z'\cup Y') \\
	x & {\rm if} \ x \in Z'\cup Y'
\end{array}
\right.\]
Although $f$ is not eventually continuous, it is sufficient to induce a functor $f_{\sharp}$ from the category $\mathcal{B}(X',Y';R)_{Z' \cap Y'}$ to the category $\mathcal{B}(Z',Z' \cap Y';R)$, defined as in Proposition $\ref{eventually}$. Furthermore, Proposition $\ref{equal}$ will apply, ensuring that $f$ is an inverse of the inclusion functor up to natural equivalence.

	We need to show that, given an object $A$ and a morphism $\phi$ in $\mathcal{B}(X',Y';R)_{Z' \cap Y'}$, $(f_{\sharp}A)_z$ is a finitely generated $R$-module, $f_{\sharp}A$ has the local finiteness property, and that $f_{\sharp}\phi$ is continuously controlled. Although $f$ is not continuous on all of $Y'$, it is continuous at $Z'\cap Y'$. Since objects are zero in a neighborhood of $Y'-Z'$, this is enough to make $f_{\sharp}\phi$ continuously controlled. (See the proof of Proposition $\ref{eventually}$.) The remaining two statements will be proved once we show that the inverse image under $f$ of a compact set in $Z'-Y'$ contains only finitely many points over which an object is non-zero.

	Given $z\in Z'-Y'$ and $0<d<d(z,Y')$, $x\in f^{-1}(\overline{B_d(z)})$ implies $d(x,Z'\cap Y')>0$. To see why, note that $d(z_x,z)\leq d$ and $d(z_x,y_x)\leq d(x,y_x)$. Therefore
\[ d(z,Y')\leq d(z,y_x)\leq d(z,z_x)+d(z_x,y_x)\leq d+d(x,y_x)=d+d(x,Z'\cap Y'). \]
Hence, $0<d(z,Y')-d\leq d(x,Z'\cap Y')$. This implies that $\overline{f^{-1}(\overline{B_d(z)})} \subseteq X'-U$, where $U=\{ x\in X' \, | \, d(x,Z'\cap Y')<d(z,Y')-d \}$. Given an object $A$ in $\mathcal{B}(X',Y';R)_{Z' \cap Y'}$, there is a neighborhood $V$ of $Y'-Z'$ on which $A$ is zero. Therefore,
\[ \big\{x\in \overline{f^{-1}(\overline{B_d(z)})} \, \big| \, A_x \neq 0 \big\} \subset (X'-U)\cap (X'-V). \]
Since $(X'-U)\cap (X'-V)$ is a compact subset of $X'-Y'$, $\big\{x\in \overline{f^{-1}(\overline{B_d(z)})} \, \big| \, A_x \neq 0 \big\}$ is finite. Since $f|_{Y'}=1_{Y'}$, the prrof is complete.
\end{proof}

\begin{lem}\label{gammaretract}
	Let $G$ be a finite group acting on a compact metrizable space $X'$, let $Y'$ be a closed, $G$-invariant, nowhere dense subset of $X'$, and let $Z'$ be a closed, $G$-invariant subset of $X'$ such that $\overline{Z'-Y'}=Z'$. Then $\mathcal{B}(X',Y';R)_{Z' \cap Y'}^G \cong \mathcal{B}(Z',Z' \cap Y';R)^G$.
\end{lem}

\begin{proof}
	We need to generalize the function $f$ described in the proof of Lemma $\ref{retract}$. Choose $x_i \in X'-(Y'-Z')$ so that $\bigsqcup_{}^{}Gx_i=X'-(Y'-Z')$. For each $i$, choose an orbit $Gy_i$ in $Z'\cap Y'$ such that $d(Gx_i,Gy_i)=d(Gx_i,Z'\cap Y')$. Since $Z'-Y'$ is dense in $Z'$, choose an orbit $Gz_i$ in $Z'-Y'$ such that $d(Gz_i,Gy_i) \leq d(Gx_i,Gy_i)$. Define $f'(x_i)=z_i$, and $f'(x)=x$ if $x\in Y'\cup Z'$. (Note that $f'$ is not defined on all of $X'$.) Recall that an object $A$ in the the fixed category has the property that $A_x\cong A_{gx}$ as $R$-modules for every $g\in G$, and that $A_{Gx}=\bigoplus_{x'\in Gx}^{}A_{x'}$ is an $RG$-module. So define the object $f'_{\sharp}A$ as follows. $(f'_{\sharp}A)_z$ is the free $RG_z$-module of rank equal to $\sum_{x \in (f')^{-1}(z)}^{}{\rm rank}_{RG}(A_{Gx})$. As with eventually continuous maps, the morphism assignments here are induced by the identity. This gives well-defined morphisms since morphisms in the fixed category are equivariant. 

	Notice that if we restrict $f$ in the proof of Lemma $\ref{retract}$ to $\{ x_i \}\cup Y'\cup Z'$, we get $f'$. Since $G$ is finite, we can assume without loss of generality that it acts on $X'$ by isometries. This implies $d(z,Y')=d(gz,Y')$ for every $g\in G$, since $Y'$ is a $G$-invariant subspace. Therefore the proof of Lemma $\ref{retract}$ tells us that $f'_{\sharp}A$ has the finiteness properties required of objects. Furthermore, as in Lemma $\ref{retract}$, morphisms are continuously controlled by a proof similar to the one given for the case of eventually continuous maps in Proposition $\ref{eventually}$.
\end{proof}

\begin{lem}\label{p_x}
	Let $G$ be a finite group acting on a compact metrizable space $X$, let $Y$ be a closed, $G$-invariant, nowhere dense subset of $X$, and let $Z$ be a closed, $G$-invariant subset of $X$ such that $\overline{Z-Y}=Z$. Then
\begin{enumerate}
\item[1.] $\mathcal{B}(CX,CY \cup X,p_X;R)_{C(Z \cap Y) \cup Z}^G \cong \mathcal{B}(CZ,C(Z \cap Y) \cup Z,p_Z;R)^G$;
\item[2.] $\mathcal{B}(\Sigma X,\Sigma Y,p_X;R)_{\Sigma Z}^G \cong \mathcal{B}(\Sigma Z,\Sigma (Z \cap Y),p_Z;R)^G$.
\end{enumerate}
\end{lem}

\begin{proof}
	Notice that the first statement is essentially the same as Lemma $\ref{gammaretract}$, taking $X'=CX$, $Y'=CY\cup X$, and $Z'=CZ$. The only difference is that we are working with morphisms that are $p_X$-controlled along $Y\times (0,1)$ instead of ones that are continuously controlled along $Y\times (0,1)$. Therefore we need only show that a morphism in the image of the functor described in Lemma~\ref{gammaretract} is $p_Z$-controlled, where $p_Z=p_X|_{(Z\cap Y) \times (0,1)}$. Specifically, we need to show that given a point $(y,t) \in (Z\cap Y) \times (0,1)$ and a neighborhood $U \subseteq Z$ of $y$, there is a neighborhood $V \subseteq CZ$ of $(y,t)$, such that $(f'_{\sharp}\phi)^a_b=0$ and $(f'_{\sharp}\phi)^b_a=0$ if $a\in V$ and $b \notin U \times (0,1)$. It suffices to prove it when $U=B_{\epsilon}(y,t)\subseteq Z$ for all $\epsilon >0$. As in the proof of Proposition $\ref{eventually}$, we do this by contradiction.

	Suppose that there is an $\epsilon >0$ and sequences $(a'_n),(b'_n)$, such that for each $n$, $a'_n \notin B_{\epsilon}(y,t) \times (0,1)$, $(b'_n)$ converges to $(y,t)$ and $(f'_{\sharp}\phi)^{a'_n}_{b'_n} \neq 0$. For each $n$, there exists $a_n$ and $b_n$ such that the $R\Gamma$-module over $Ga_n$ is sent to $Ga'_n$, the $R\Gamma$-module over $Gb_n$ is sent to $Gb'_n$, and $\phi^{a_n}_{b_n} \neq 0$. Since $X$ is compact Hausdorff, $(b_n)$ has a convergent subsequence $(b_k)$ that converges to some point $b$. The local finiteness condition on objects implies $b\in C(Z\cap Y)\cup Z$. By the definition of our functor, this implies that $d(Gb'_n,Gb_n)$ converges to zero. Hence $b=(y,t)$. The definition of the functor also implies $d(a_k,(y,t))\geq \epsilon$ for every $k$. This means that $a_k$ is outside of $B_{\epsilon}(y,t) \times (0,1) \subseteq CX$ for every $k$. But that contradicts the $p_X$-control of $\phi$. This completes the proof of the first part of this lemma. The second part is proved precisely the same way, replacing cones with suspensions.
\end{proof}

\section{Homology and Continuously Controlled Algebra}\label{hcca}

	In \cite{karoubi}, Karoubi introduced the notion of an $\mathcal{A}$-{\it filtered} additive category $\mathcal{U}$.

\begin{defi}
	{\rm Let $\mathcal{A}$ be a full subcategory of an additive category $\mathcal{U}$. Denote the objects of $\mathcal{A}$ by the letters $A$ through $F$ and the objects of $\mathcal{U}$ by the letters $U$ through $W$. Then $\mathcal{U}$ is said to be $\mathcal{A}$-{\it filtered} if every object $U$ has a family of decompositions $\{ U=E_{\alpha} \oplus U_{\alpha} \}$, where
\begin{enumerate}
	\item[(i)] the decomposition forms a filtered poset under the partial order in which $E_{\alpha} \oplus U_{\alpha} \leq E_{\beta} \oplus U_{\beta}$ whenever $U_{\beta} \subseteq U_{\alpha}$ and $E_{\alpha} \subseteq E_{\beta}$;
	\item[(ii)] every map $A \to U$ factors through $E_{\alpha}$ for some $\alpha$;
	\item[(iii)] every map $U \to A$ factors through $E_{\alpha}$ for some $\alpha$;
	\item[(iv)] for each $U$ and $V$, the filtration on $U \oplus V$ is equivalent to the sum of the filtrations $\{ U=E_{\alpha} \oplus U_{\alpha} \}$ and $\{ V=F_{\beta} \oplus V_{\beta} \}$; that is, $\{ U \oplus V = (E_{\alpha} \oplus F_{\beta}) \oplus (U_{\alpha} \oplus V_{\beta}) \}$.
\end{enumerate} }
\end{defi}

	Karoubi also defined $\mathcal{U} / \mathcal{A}$ to be the category with the same objects as $\mathcal{U}$, but morphisms $\phi$, $\psi : U \to V$ are identified if $\phi - \psi$ factors through $\mathcal{A}$.

	Let $\mathcal{A}$ be a small additive category.  Let ${\mathbb{K}}^{-\infty}(\mathcal{A})$ denote the non-connective $K$-theory spectrum associated with the symmetric monoidal category obtained from $\mathcal{A}$ by restricting to isomorphisms. Then ${\mathbb{K}}^{-\infty}$ is a functor from the category of small additive categories to the category of spectra. C\'ardenas and Pedersen give a detailed description of this functor in \cite{cardenas}.

	Karoubi filtrations become the main tool in the theory of continuously controlled algebra via Theorem~\ref{fibration} below. A proof of it can also be found in \cite{cardenas}.

\begin{theo}\cite[Theorem 1.28]{cp}\label{fibration}
	The sequence
\[ {\mathbb{K}}^{-\infty}(\mathcal{A}) \to {\mathbb{K}}^{-\infty}(\mathcal{U}) \to {\mathbb{K}}^{-\infty}(\mathcal{U} / \mathcal{A}) \]
is a homotopy fibration of spectra.
\end{theo}

\begin{lem}\label{equivariantK}
	Let $\Gamma$ be a group acting on a compact Hausdorff space $X$, let $Y$ be a closed $\Gamma$-invariant subspace such that $X-Y$ is dense in $X$, let $C$ be a closed $\Gamma$-invariant subspace of $Y$, and let $W$ an open $\Gamma$-invariant subset of $Y$ containing $C$. Then 
\begin{enumerate}
	\item[1.] $\mathcal{B}(X,Y;R)^{\Gamma}$ is $\mathcal{B}(X,Y;R)_C^{\Gamma}$-filtered;
	\item[2.] $\mathcal{B}(X,Y;R)^{\Gamma} / \mathcal{B}(X,Y;R)_C^{\Gamma} = (\mathcal{B}(X,Y;R)^{Y-C})^{\Gamma}$;
	\item[3.] $(\mathcal{B}(X,Y;R)^W)^{\Gamma}$ is $(\mathcal{B}(X,Y;R)_C^W)^{\Gamma}$-filtered;
	\item[4.] $(\mathcal{B}(X,Y;R)^W)^{\Gamma} / (\mathcal{B}(X,Y;R)_C^W)^{\Gamma} = (\mathcal{B}(X,Y;R)^{W-C})^{\Gamma}$.
\end{enumerate}
\end{lem}

\begin{proof}
	This follows immediately from the definitions.
\end{proof}

\begin{coro}
	Under the conditions of Lemma $\ref{equivariantK}$, the following sequences are fibrations of spectra up to homotopy:
	\begin{enumerate}
	\item[1.] ${\mathbb{K}}^{-\infty}(\mathcal{B}(X,Y;R)_C)^{\Gamma} \to {\mathbb{K}}^{-\infty}(\mathcal{B}(X,Y;R))^{\Gamma} \to {\mathbb{K}}^{-\infty}(\mathcal{B}(X,Y;R)^{Y-C})^{\Gamma}$
	\item[2.] ${\mathbb{K}}^{-\infty}(\mathcal{B}(X,Y;R)_C^W)^{\Gamma} \to {\mathbb{K}}^{-\infty}(\mathcal{B}(X,Y;R)^W)^{\Gamma} \to {\mathbb{K}}^{-\infty}(\mathcal{B}(X,Y;R)^{W-C})^{\Gamma}$.
	\end{enumerate}
\end{coro}

\begin{proof}
	This follows from Lemma $\ref{equivariantK}$ and Theorem $\ref{fibration}$ since taking fixed sets commutes with applying ${\mathbb{K}}^{-\infty}$.
\end{proof}

\begin{defi}
	{\rm A {\it reduced Steenrod homology theory}, $h$, is a functor from the category of compact metrizable spaces and continuous maps to the category of graded abelian groups satisfying
\begin{enumerate}
	\item[(i)] $h$ is homotopy invariant;
	\item[(ii)] $h(\bullet)=0$;
	\item[(iii)] given any closed subset $A$ of $X$, there is a natural transformation, $\partial :h_n(X / A) \to h_{n-1}(A)$, fitting into a long exact sequence
\[ \dotsb \to h_n(A) \to h_n(X) \to h_n(X / A) \to h_{n-1}(A) \to \dotsb ; \]
	\item[(iv)] if $\bigvee^{}_{} X_i$ denotes a compact metric space that is a countable union of metric spaces along a single point, then the projection maps $p_i : \bigvee^{}_{} X_i \to X_i$ induce an isomorphism $h_*(\bigvee^{}_{} X_i) \to \prod^{}_{} h_*(X_i)$.
\end{enumerate} }
\end{defi}

	These axioms are often called the {\it Kaminker-Schochet axioms}.

	Given any generalized homology theory, there is a unique Steenrod homology extension. Existence of such extensions was proved by Kahn, Kaminker, and Schochet and Edwards and Hastings~\cite{kahn,edwards}. Uniqueness was proved by Milnor \cite{milnor}.

\begin{defi}
	{\rm A functor $k$ from the category of compact metrizable spaces to the category of spectra is called a {\it reduced Steenrod functor} if it satisfies the following conditions.
\begin{enumerate}
	\item[(i)] The spectrum $k(CX)$ is contractible.
	\item[(ii)] If $A \subset X$ is closed, then $k(A) \to k(X) \to k(X / A)$ is a fibration (up to natural weak homotopy equivalence).
	\item[(iii)] The projection maps $p_i : \bigvee^{}_{} X_i \to X_i$ induce a weak homotopy equivalence $k(\bigvee^{}_{} X_i) \to \prod^{}_{} k(X_i)$.
\end{enumerate} }
\end{defi}

\begin{theo}\cite[Proposition 1.35]{cp}
	Let $k$ be a reduced Steenrod functor. Then ${\pi}_*(k(X))$ is the unique reduced Steenrod homology theory associated with the spectrum $k(S^0)$. If $X$ is a finite CW-complex, then $k(X)$ is weakly homotopy equivalent to $X \wedge k(S^0)$.
\end{theo}

\begin{theo}\cite[Theorem 1.36]{cp}\label{sthom}
        The functor ${\mathbb{K}}^{-\infty}(\mathcal{B}(C(-),-;R))$ is a reduced Steenrod functor whose value on $S^0$ is $\Sigma {\mathbb{K}}^{-\infty}(R)$. Furthermore, if $X$ is a finite CW-complex then $\Omega {\mathbb{K}}^{-\infty}(\mathcal{B}(CX,X;R))$ is weakly homotopy equivalent to $X \wedge {\mathbb{K}}^{-\infty}(R)$.
\end{theo}

        Clever use of Theorem~\ref{fibration} yields the following result \cite{anderson, hp}.

\begin{theo}\label{homology}
	The functor ${\mathbb{K}}^{-\infty}\left(\mathcal{B}(- \times (0,1],-\times 1;R)^{-\times 1}\right)^{\Gamma}$, from the category of $\Gamma$-spaces to the category of spectra, is $\Gamma$-homotopy invariant and $\Gamma$-excisive.
\end{theo}

	Theorem $\ref{homology}$ implies that $\pi_*({\mathbb{K}}^{-\infty}(\mathcal{B}(- \times (0,1],-\times 1;R)^{-\times 1})^{\Gamma})$ satisfies the Eilenberg-Steenrod axioms, except for the dimension axiom, and is therefore a generalized homology theory. We will denote ${\pi}_*(\Omega{\mathbb{K}}^{-\infty}(\mathcal{B}(E \times (0,1],E\times 1;R)^{E\times 1})^{\Gamma})$ by $H_*^{\Gamma} (E;{\mathbb{K}}^{-\infty}(R{\Gamma}_x))$.

	We conclude this section by connecting the category $\mathcal{B} (CX,CY \cup X,p_X;R)$ (used in the proof of Theorem~\ref{main}) to Theorem~\ref{homology}.

\begin{lem}\cite[Lemma 2.4]{cp}\label{lem3a}
	$\mathcal{B} (CX,CY \cup X,p_X;R)^E = \mathcal{B} (E \times (0,1],E \times 1;R)^{E \times 1}$.
\end{lem}

\begin{lem}\cite[Lemma 2.5]{cp}\label{lem3b}
	Let $\Gamma$ be a group acting on a compact Hausdorff space $X$, and let $Y$ be a closed $\Gamma$-invariant subspace such that $X-Y$ is dense in $X$. Then ${\mathbb{K}}^{-\infty}(\mathcal{B} (CX,CY \cup X,p_X;R))^{\Gamma} \simeq {\mathbb{K}}^{-\infty}(\mathcal{B}(CX,CY \cup X,p_X;R)^E)^{\Gamma}$.
\end{lem}

\begin{coro}\label{isom}
	Let $\Gamma$ be a group acting on a compact Hausdorff space $X$, and let $Y$ be a closed $\Gamma$-invariant subspace such that $X-Y$ is dense in $X$. Then ${\mathbb{K}}^{-\infty}(\mathcal{B} (CX,CY \cup X,p_X;R))^{\Gamma} \simeq {\mathbb{K}}^{-\infty}(\mathcal{B} (E \times (0,1],E \times 1;R)^{E \times 1})^{\Gamma}$.
\end{coro}

\begin{proof}
	This follows immediately from Lemmas $\ref{lem3a}$ and $\ref{lem3b}$.
\end{proof}

         Corollary $\ref{isom}$ implies that
\[ {\pi}_*({\mathbb{K}}^{-\infty}(\mathcal{B} (CX,CY \cup X,p_X;R))^{\Gamma}) \cong H_*^{\Gamma} (E;{\mathbb{K}}^{-\infty}(R{\Gamma}_x)). \]

\section{Generalized Fixed Sets}\label{gfs}

	We begin this section with some facts about generalized homotopy fixed sets in the category of spectra. Let $S$ be a spectrum with $\Gamma$-action. The fixed set $S^{\Gamma}$ can be identified with the set $\Map_{\Gamma}(\bullet,S)$ of $\Gamma$-equivariant maps from a point into $S$. Similarly, the homotopy fixed set $S^{h\Gamma}$ is defined to be $\Map_{\Gamma}(E\Gamma,S)$, where $E\Gamma$ denotes the universal space for free $\Gamma$-actions. It is a well-known theorem that an equivariant map, $S\to T$, that is a weak homotopy equivalence unequivariantly, induces a weak homotopy equivalence on homotopy fixed sets if $E\Gamma$ is assumed to be a finite $\Gamma$-CW complex.

	Let $\mathcal{F}$ be a family of subgroups of $\Gamma$ that is closed under conjugation and under the operation of taking subgroups. Then there is a universal space $E\Gamma(\mathcal{F})$ for $\Gamma$-actions with isotropy in $\mathcal{F}$. The space $E\Gamma(\mathcal{F})$ is a $\Gamma$-CW complex that is characterized by the fact that $E\Gamma(\mathcal{F})^G$ is contractible for all $G\in \mathcal{F}$ and is empty otherwise. Denote $\Map_{\Gamma}(E\Gamma(\mathcal{F}),S)$ by $S^{h \mathcal{F} \Gamma}$. The next two lemmas play an important role in the proof of Theorem $\ref{main}$. Their proofs come from unpublished work of Carlsson, Pedersen and Roe.

\begin{lem}\label{lemmaa}
	Let $F:S \to T$ be an equivariant map of spectra with $\Gamma$-action. Assume that $E\Gamma(\mathcal{F})$ is a finite $\Gamma$-CW complex. If $F^G:S^G \to T^G$ is a weak homotopy equivalence for every $G\in \mathcal{F}$, then $S^{h \mathcal{F} \Gamma} \simeq T^{h \mathcal{F} \Gamma}$.
\end{lem}

\begin{proof}
	Proceed by induction on the $\Gamma$-cells in $E\Gamma(\mathcal{F})$. Consider the $\Gamma$-$0$-cell, $\Gamma / H$, where $H\in \mathcal{F}$. Then $\Map_{\Gamma}(\Gamma / H,S)=S^H \simeq T^H=\Map_{\Gamma}(\Gamma / H,T)$, which takes care of the base case.

	Assume now that $X$ is obtained from $A$ by attaching a $\Gamma$-$n$-cell, $\Gamma / K \times \mathbb{D}^n$, for some $K\in \mathcal{F}$, and that $\Map_{\Gamma}(A,S) \simeq \Map_{\Gamma}(A,T)$. We need to show that $\Map_{\Gamma}(X,S) \simeq \Map_{\Gamma}(X,T)$. Note that $A \to X \to X / A$ is a $\Gamma$-equivariant cofibration sequence. This implies that $\Map_{\Gamma}(X / A,S) \to \Map_{\Gamma}(X,S) \to \Map_{\Gamma}(A,S)$ is a fibration of spectra. By the induction hypothesis, $\Map_{\Gamma}(A,S) \simeq \Map_{\Gamma}(A,T)$. Therefore, once we show that $\Map_{\Gamma}(X / A,S) \simeq \Map_{\Gamma}(X / A,T)$, the Five Lemma will finish the proof.

	Now, since $X / A = {\Sigma}^n (\Gamma / K)$, by adjointness it suffices to show that
\[ \Map_{\Gamma}(\Gamma / K,{\Omega}^n S) \simeq \Map_{\Gamma}(\Gamma / K,{\Omega}^n T). \]
This is equivalent to proving $({\Omega}^n S)^K \simeq ({\Omega}^n T)^K$. But this holds since $({\Omega}^n S)^K = {\Omega}^n (S^K)$ and $S^K \simeq T^K$.
\end{proof}

\begin{lem}\label{lemmab}
	Let $B$ be a $G$-spectrum, where $G\in \mathcal{F}$, and let $\Gamma$ act on $S=\prod_{\Gamma / G}B$ by identifying $\prod_{\Gamma / G}B$ with $\Map_G(\Gamma,B)$, in which $(\gamma f)(x)=f({\gamma}^{-1}x)$, where $\gamma \in \Gamma$. Then $S^{\Gamma} \simeq S^{h \mathcal{F} \Gamma}$.
\end{lem}

\begin{proof}
\[ S^{h \mathcal{F} \Gamma}= \Map_{\Gamma}(E\Gamma(\mathcal{F}),\Map_G(\Gamma,B))=\Map_G(E\Gamma(\mathcal{F}),B) \]
by evaluating at 1.
Since $E\Gamma(\mathcal{F})$ is $G$-equivariantly homotopy equivalent to $E\Gamma(\mathcal{F})^G$ and $E\Gamma(\mathcal{F})^G$ is contractible,
\begin{align}
	&\Map_G(E\Gamma(\mathcal{F}),B) & &\cong & &\Map_G(E\Gamma(\mathcal{F})^G,B) & &= & &\Map(E\Gamma(\mathcal{F})^G,B^G) \notag\\
	& & & & & & &\cong & &B^G\notag\\
	& & & & & & &= & &S^{\Gamma}.\notag
\end{align}
\end{proof}

	Notice that since we can identify $S^{\Gamma}$ with $\Map_{\Gamma}(\bullet,S)$, the $\Gamma$-equivariant map $E\Gamma(\mathcal{F}) \to \bullet$
induces the following commutative diagram:
\[ \xymatrix{
	S^{\Gamma} \ar[d] \ar[r] & T^{\Gamma} \ar[d] \\
	S^{h\mathcal{F} \Gamma} \ar[r] & T^{h\mathcal{F} \Gamma} } \]
Such a diagram will be used in the proof of Theorem~\ref{main}.

\section{The Conner Conjecture and the Vanishing of Steenrod Homology}\label{ccvsh}	

	In 1960 Conner conjectured that the orbit space of any action of a compact Lie group on Euclidean $n$-space~--- or on the closed $n$-disk~--- is contractible \cite{conner}. This conjecture was motivated by a result of Floyd in 1951 that whenever a finite group operates as a group of topological transformations on a compact finite dimensional ANR, then the orbit space is also an ANR~\cite{floyd}. To prove his conjecture, Conner further conjectured that if a compact Lie group $G$ acts on a space $X$, where $X$ is either paracompact of finite cohomological dimension with finitely many orbit types or compact Hausdorff, then $X / G$ has trivial reduced \v{C}ech cohomology if $X$ does \cite{conner}. This latter conjecture was proved by Oliver in 1976~\cite{oliver}.

	In Section~\ref{fccfs}, we will want to know that the reduced Steenrod homology (from Theorem~\ref{sthom}) of the quotient of a contractible compact metrizable space by a finite group is trivial. The first step toward proving this is Theorem~\ref{conconj} below, which is a previously known special case of Oliver's result. A proof can be found in~\cite{bredon}.

\begin{theo}\cite[Theorem III.7.12]{bredon}\label{conconj}
	Let $G$ be a finite group acting on a compact Hausdorff space $X$. If the reduced \v{C}ech cohomology of $X$ is trivial, then the reduced \v{C}ech cohomology of $X / G$ is also trivial.
\end{theo}

	The final step is the following theorem.

\begin{theo}\label{geogh}
	If the reduced \v{C}ech cohomology of a compact metrizable space $X$ is trivial, then the reduced Steenrod homology (from Theorem $\ref{sthom}$) of $X$ is also trivial.
\end{theo}

\begin{proof}
	Let $\Gamma$ be an abelian group, let $H_q(X;\Gamma)$ denote the ordinary unreduced Steenrod homology of $X$ with coefficients in $\Gamma$, let $\widetilde{H}_q(X;\Gamma)$ denote the ordinary reduced Steenrod homology of $X$ with coefficients in $\Gamma$, and let $\check{H}^q(X;\Gamma)$ denote the unreduced \v{C}ech cohomology of $X$ with coefficients in $\Gamma$. Note that $H$ is related to $\check{H}$ by a split exact sequence
\[ 0 \to {\rm Ext}(\check{H}^{q+1}(X),\Gamma) \to H_q(X;\Gamma) \to {\rm Hom}(\check{H}^q(X),\Gamma) \to 0 \]
 \cite[Lemma 5]{milnor}. Since the reduced \v{C}ech cohomology of $X$ is trivial, this implies that $\widetilde{H}_q(X;\Gamma)=0$ for all $q$.

	There is a spectral sequence converging to ${\pi}_*({\mathbb{K}}^{-\infty}(\mathcal{B}(CX,X;R)))$ with
\[ E^2_{p,q}=\widetilde{H}_p\left(X;{\pi}_q({\mathbb{K}}^{-\infty}(\mathcal{B}(C(S^0),S^0;R)))\right)=\widetilde{H}_p\left(X;K_{q-1}(R)\right) \]
\cite[8.5.5]{edwards}. Since each term of the spectral sequence is zero, the proof is complete.
\end{proof}

\begin{coro}\label{zero}
	Let $G$ be a finite group acting on a compact metrizable space $X$. If the reduced \v{C}ech cohomology of $X$ is trivial, then the reduced Steenrod homology of $X / G$ is trivial.
\end{coro}

	The reason for wanting a space to have trivial reduced Steenrod homology is so that we may use the following proposition. It follows from \cite[Proof of Theorem 2.13]{cp}.

\begin{prop}\label{key}
	Let $X$ be a compact metrizable space and $Y$ a closed nowhere dense subset. If the reduced Steenrod homology of $X$ is trivial, then
\[ {\mathbb{K}}^{- \infty} (\mathcal{B} (CX,CY \cup X,p_X;R)) \simeq {\mathbb{K}}^{- \infty} (\mathcal{B} (\Sigma X,\Sigma Y,p_X;R)). \]
\end{prop}

\section{Algebraic $K$-theory and Continuously Controlled Algebra}\label{akcca}

        In this section we state the main result of this paper (Theorem~\ref{main} below) and begin its proof.

\begin{theo}\label{main}
	Let $\Gamma$ be a discrete group, and let $\mathcal{E}=E\Gamma(\mathfrak{f})$ be the universal space for $\Gamma$-actions with finite isotropy, where $\mathfrak{f}$ denotes the family of finite subgroups of $\Gamma$. Assume that $\mathcal{E}$ is a finite $\Gamma$-CW complex admitting a compactification, $X$, (i.e., $X$ is compact and $\mathcal{E}$ is an open dense subset) such that
	\begin{enumerate}
		\item[1.] the $\Gamma$-action extends to $X$;
		\item[2.] $X$ is metrizable;
		\item[3.] $X^G$ is contractible for every $G\in \mathfrak{f}$;
		\item[4.] ${\mathcal{E}}^G$ is dense in $X^G$ for every $G\in \mathfrak{f}$;
		\item[5.] compact subsets of $\mathcal{E}$ become small near $Y=X-\mathcal{E}$.  That is, for every compact subset $K \subset \mathcal{E}$ and for every neighborhood $U \subset X$ of $y \in Y$, there exists a neighborhood $V \subset X$ of $y$ such that $g \in \Gamma$ and $gK \cap V \neq \emptyset$ implies $gK \subset U$.
	\end{enumerate}
Then $H_*^{\Gamma} (\mathcal{E};{\mathbb{K}}^{-\infty}(R{\Gamma}_x)) \to K_*(R\Gamma)$ is split injective for every ring $R$.
\end{theo}

\begin{exam}
	{\rm A {\it crystallographic group} is a discrete group that acts cocompactly on Euclidean $n$-space by isometries. Let $\Gamma$ be a crystallographic group acting on $\mathbb{R}^n$. Then $\mathbb{R}^n$ is a universal space for $\Gamma$-actions with finite isotropy. Furthermore, $\mathbb{R}^n$ is a finite $\Gamma$-CW complex whose fixed sets are totally geodesic. This implies that if we compactify $\mathbb{R}^n$ with an $(n-1)$-sphere at infinity, then $\Gamma$ will satisfy all of the conditions of Theorem~\ref{main}. More generally, virtually polycyclic groups satisfy all of the conditions of Theorem~\ref{main}, since we can also take some $\mathbb{R}^n$ to be our universal space. This is discussed in greater detail in Section~\ref{vpg}.} 
\end{exam}

	For the remainder of this section and throughout Section~\ref{fccfs}, assume that $\Gamma$, $\mathcal{E}$, $X$, and $Y$ satisfy the hypotheses of Theorem~\ref{main}.
	
	Consider the $\Gamma$-equivariant functor \[\mathcal{B} (CX,CY \cup X,p_X;R) \to \mathcal{B} (\Sigma X,\Sigma Y,p_X;R) \] induced by collapsing $X$ from $CX$ to $\Sigma X$. Let $S = \Omega {\mathbb{K}}^{- \infty} (\mathcal{B} (CX,CY \cup X,p_X;R))$ and $T = \Omega {\mathbb{K}}^{- \infty} (\mathcal{B} (\Sigma X,\Sigma Y,p_X;R))$. Recall that ${\pi}_*(S^{\Gamma}) \cong H_*^{\Gamma} (\mathcal{E};{\mathbb{K}}^{-\infty}(R{\Gamma}_x))$ by Theorem $\ref{homology}$ and Corollary $\ref{isom}$. It is proved in \cite[Lemma 2.3]{cp} that $T^{\Gamma} \simeq {\mathbb{K}}^{- \infty}(R \Gamma)$. Thus, the map of fixed sets, $S^{\Gamma} \to T^{\Gamma}$, is the object of focus. This formulation of the continuously controlled assembly map differs from the one presented in~\cite{hp}, however the two are homotopy equivalent. The benefit of using the map described here is that it is a map a fixed sets which allows us to use the following commutative diagram:
\[ \xymatrix{
	S^{\Gamma} \ar[d] \ar[r] & T^{\Gamma} \ar[d] \\
	S^{h\mathfrak{f} \Gamma} \ar[r] & T^{h\mathfrak{f} \Gamma}. } \]
Then Theorem $\ref{main}$ is proved by the following two statements:

\begin{enumerate}
	\item[1.] $S^{\Gamma} \simeq S^{h \mathfrak{f} \Gamma}$;
	\item[2.] $S^{h \mathfrak{f} \Gamma} \simeq T^{h \mathfrak{f} \Gamma}$.
\end{enumerate}
	To prove the first statement, we need a slight generalization of the results in \cite[Theorem 2.11, Corollary 2.12]{cp}. 

\begin{theo}\label{theo3}
	The spectrum ${\mathbb{K}}^{-\infty}\left(\mathcal{B} (\mathcal{E} \times (0,1],\mathcal{E} \times 1;R)^{\mathcal{E} \times 1}\right)^{\Gamma}$ is weakly homotopy equivalent to ${\mathbb{K}}^{-\infty}\left(\mathcal{B} (\mathcal{E} \times (0,1],\mathcal{E} \times 1;R)^{\mathcal{E} \times 1}\right)^{h \mathfrak{f} \Gamma}$.
\end{theo}

\begin{proof}
	Note that $\mathcal{E}$ is assumed to be a finite $\Gamma$-CW complex, and proceed by induction on the $\Gamma$-cells in $\mathcal{E}$.  Begin with the discrete space ${\Gamma}/H$ for some $H \in \mathfrak{f}$. The control condition implies that the components of a morphism near ${\Gamma}/H \times 1$ must be zero between points with different ${\Gamma}/H$ entries. Since we are taking germs at ${\Gamma}/H \times 1$, the category $\mathcal{B} ({\Gamma}/H \times (0,1],{\Gamma}/H \times 1;R)^{{\Gamma}/H \times 1}$ is equivalent to the product category $\prod_{\Gamma /H}^{} \mathcal{B} ((0,1],1;R)^{1}$. The action of $\Gamma$ that the product category inherits is the same as the one on the product in Lemma~\ref{lemmab}. The projection maps induce a map
\[ {\mathbb{K}}^{-\infty} \Big(\prod_{\Gamma /H}^{} \mathcal{B} ((0,1],1;R)^{1}\Big) \to \prod_{\Gamma /H}^{} {\mathbb{K}}^{-\infty} \left(\mathcal{B} ((0,1],1;R)^{1}\right) \]
that is $\Gamma$-equivariant, and a weak homotopy equivalence by~\cite{carlsson}.

	Consider the following commutative diagram:
\[ \xymatrix{
	{{\mathbb{K}}^{-\infty} \left(\prod_{\Gamma /H}^{} \mathcal{B} ((0,1],1;R)^{1}\right)}^{\Gamma} \ar[d]_a \ar[r]^b &  {\left( \prod_{\Gamma /H}^{} {\mathbb{K}}^{-\infty} \left(\mathcal{B} ((0,1],1;R)^{1}\right) \right)}^{\Gamma} \ar[d]_c \\
	{{\mathbb{K}}^{-\infty} \left(\prod_{\Gamma /H}^{} \mathcal{B} ((0,1],1;R)^{1}\right)}^{h\mathfrak{f} \Gamma} \ar[r]^d & {\left( \prod_{\Gamma /H}^{} {\mathbb{K}}^{-\infty} \left(\mathcal{B} ((0,1],1;R)^{1}\right) \right)}^{h\mathfrak{f} \Gamma} } \]
We want to show that $a$ is a weak homotopy equivalence.

	Taking fixed sets commutes with applying ${\mathbb{K}}^{-\infty}$, therefore
\begin{align}
	&{\mathbb{K}}^{-\infty} \Big(\prod_{\Gamma /H}^{} \mathcal{B} ((0,1],1;R)^{1}\Big)^{\Gamma} & &\cong & &{\mathbb{K}}^{-\infty} \Big( \Big( \prod_{\Gamma /H}^{} \mathcal{B} ((0,1],1;R)^{1} \Big)^{\Gamma} \Big) \notag\\ 
	& & &\cong & &{\mathbb{K}}^{-\infty} \left( \big(\mathcal{B} ((0,1],1;RH)^{1}\big) \right) \notag\\
	& & &\cong & &{\mathbb{K}}^{-\infty} \left( \big(\mathcal{B} ((0,1],1;R)^{1}\big)^H \right) \notag\\
	& & &\cong & &{\mathbb{K}}^{-\infty} \left( \mathcal{B} ((0,1],1;R)^{1} \right)^H \notag\\
	& & &\cong & &\Big( \prod_{\Gamma /H}^{} {\mathbb{K}}^{-\infty} \left(\mathcal{B} ((0,1],1;R)^{1}\right) \Big)^{\Gamma}. \notag
\end{align}
That is, $b$ is a homotopy equivalence. More generally, if $G \in \mathfrak{f}$, choose representatives, $\gamma_i$ so that $G \backslash \Gamma /H=\{ G\gamma_iH \}$. Then 
\begin{align}
	&{\mathbb{K}}^{-\infty} \Big(\prod_{\Gamma /H} \mathcal{B} ((0,1],1;R)^{1}\Big)^G & &\cong & &{\mathbb{K}}^{-\infty} \Big( \Big( \prod_{\Gamma /H}^{} \mathcal{B} ((0,1],1;R)^{1} \Big)^G \Big) \notag\\ 
	& & &\cong & &{\mathbb{K}}^{-\infty} \Big( \prod_{G \backslash \Gamma /H} \mathcal{B} ((0,1],1;R[\gamma_i^{-1}G\gamma_i \cap H])^{1} \Big) \notag\\
	& & &\simeq & &\prod_{G \backslash \Gamma /H} {\mathbb{K}}^{-\infty} \left( \mathcal{B} ((0,1],1;R)^{1} \right)^{\gamma_i^{-1}G\gamma_i \cap H} \notag\\
	& & &\cong & &\Big( \prod_{\Gamma /H} {\mathbb{K}}^{-\infty} \left(\mathcal{B} ((0,1],1;R)^{1}\right) \Big)^G \notag
\end{align}
since $\mathbb{K}^{-\infty}$ commutes with taking infinite products up to weak homotopy equivalence~\cite{carlsson}.
This implies that $d$ is a weak homotopy equivalence by Lemma $\ref{lemmaa}$. Finally, $c$ is a weak homotopy equivalence by Lemma $\ref{lemmab}$.

	Now assume that the theorem holds for $N$ and that $E$ is obtained from $N$ by attaching a $\Gamma$-$n$-cell, $\Gamma / K \times \mathbb{D}^n$, for some $K \in \mathfrak{f}$. Since $\mathcal{B}(E \times (0,1],E \times 1;R)_{N \times 1}^{E \times 1}$ is $\Gamma$-equivariantly equivalent to $\mathcal{B}(N \times (0,1],N \times 1;R)^{N \times 1}$, $\mathcal{B}(N \times (0,1],N \times 1;R)^{N \times 1} \to \mathcal{B}(E \times (0,1],E \times 1;R)^{E \times 1} \to \mathcal{B}(E \times (0,1],E \times 1;R)^{(E-N) \times 1}$ is a Karoubi filtration. Let $A \to B \to C$ denote ${\mathbb{K}}^{-\infty}$ applied to this sequence. Then $A \to B \to C$ and $A^{\Gamma} \to B^{\Gamma} \to C^{\Gamma}$ are both fibrations of spectra. Since $N$ is $\Gamma$-invariant, it is also $G$-invariant for every $G \in \mathfrak{f}$. Hence, $A^G \to B^G \to C^G$ is a fibration for every $G \in \mathfrak{f}$. Let $D$ be the homotopy fiber of $B \to C$. Taking fixed sets and taking homotopy fibers commute since both are inverse limits. Therefore $D^G \to B^G \to C^G$ is a fibration for every $G \in \mathfrak{f}$. Furthermore $A^G \simeq D^G$ for every $G \in \mathfrak{f}$. So, by Lemma $\ref{lemmaa}$, $A^{h \mathfrak{f} \Gamma} \simeq D^{h \mathfrak{f} \Gamma}$. Taking homotopy fixed sets is also an inverse limit. Therefore, it too commutes with taking homotopy fibers. This implies that $D^{h \mathfrak{f} \Gamma} \to B^{h \mathfrak{f} \Gamma} \to C^{h \mathfrak{f} \Gamma}$ is a fibration, and, as a result, so is $A^{h \mathfrak{f} \Gamma} \to B^{h \mathfrak{f} \Gamma} \to C^{h \mathfrak{f} \Gamma}$.
	
	Consider the following commutative diagram:
\[ \xymatrix{
	A^{\Gamma} \ar[d]_a \ar[r] & B^{\Gamma} \ar[d]_b \ar[r] & C^{\Gamma} \ar[d]_c \\
	A^{h \mathfrak{f} \Gamma} \ar[r] & B^{h \mathfrak{f} \Gamma} \ar[r] & C^{h \mathfrak{f} \Gamma} } \]
We need to show that $b$ is a weak homotopy equivalence. By our induction hypothesis, $a$ is a weak homotopy equivalence. So, by the Five Lemma, it suffices to prove that $c$ is a weak homotopy equivalence. Since $E-N= \Gamma / K \times \mathring{e^n}$, where $\mathring{e^n}$ is an open $n$-cell, the category $\mathcal{B} (E \times (0,1],E \times 1;R)^{(E-N) \times 1}$ is equivalent to the product category $\prod_{\Gamma /K}^{} \mathcal{B} (\mathring{e^n} \times (0,1],\mathring{e^n} \times 1;R)^{\mathring{e^n} \times 1}$. But this is entirely similar to the start of the induction. This completes the proof.
\end{proof}

\begin{coro}
	$S^{\Gamma} \simeq S^{h \mathfrak{f} \Gamma}$.
\end{coro}

\begin{proof}
	Using Lemma $\ref{lem3a}$, we have ${\mathbb{K}}^{-\infty} (\mathcal{B} (CX,CY \cup X,p_X;R)^\mathcal{E})^{\Gamma} = {\mathbb{K}}^{-\infty}(\mathcal{B} (\mathcal{E} \times (0,1],\mathcal{E} \times 1;R)^{\mathcal{E} \times 1})^{\Gamma}$ and ${\mathbb{K}}^{-\infty}(\mathcal{B} (CX,CY \cup X,p_X;R)^\mathcal{E})^{h \mathfrak{f} \Gamma} = {\mathbb{K}}^{-\infty}(\mathcal{B} (\mathcal{E} \times (0,1],\mathcal{E} \times 1;R)^{\mathcal{E} \times 1})^{h \mathfrak{f} \Gamma}$. Lemmas $\ref{lemmaa}$ and $\ref{lem3b}$ imply that ${\mathbb{K}}^{-\infty}(\mathcal{B} (CX,CY \cup X,p_X;R))^{h \mathfrak{f} \Gamma} \simeq {\mathbb{K}}^{-\infty}(\mathcal{B} (CX,CY \cup X,p_X;R)^\mathcal{E})^{h \mathfrak{f} \Gamma}$. By Lemma $\ref{lem3b}$ and Theorem~\ref{theo3}, 
\begin{align}
	&{\mathbb{K}}^{-\infty}(\mathcal{B} (CX,CY \cup X,p_X;R))^{\Gamma} & &\simeq & &{\mathbb{K}}^{-\infty}\left(\mathcal{B} (CX,CY \cup X,p_X;R)^\mathcal{E}\right)^{\Gamma} \notag\\ 
	& & &= & &{\mathbb{K}}^{-\infty}\left(\mathcal{B} (\mathcal{E} \times (0,1],\mathcal{E} \times 1;R)^{\mathcal{E} \times 1}\right)^{\Gamma} \notag\\
	& & &\simeq & &{\mathbb{K}}^{-\infty}\left(\mathcal{B} (\mathcal{E} \times (0,1],\mathcal{E} \times 1;R)^{\mathcal{E} \times 1}\right)^{h \mathfrak{f} \Gamma} \notag\\
	& & &= & &{\mathbb{K}}^{-\infty}\left(\mathcal{B} (CX,CY \cup X,p_X;R)^\mathcal{E}\right)^{h \mathfrak{f} \Gamma} \notag\\
	& & &\simeq & &{\mathbb{K}}^{-\infty}(\mathcal{B} (CX,CY \cup X,p_X;R))^{h \mathfrak{f} \Gamma}. \notag
\end{align}
Therefore $S^{\Gamma} \simeq S^{h \mathfrak{f} \Gamma}$.
\end{proof}


\section{Filtering by Conjugacy Classes of Fixed Sets}\label{fccfs}

	In this section we prove that $S^{h \mathfrak{f} \Gamma} \simeq T^{h \mathfrak{f} \Gamma}$, which will complete the proof of Theorem~\ref{main}. In light of Lemma $\ref{lemmaa}$, it suffices to prove the following theorem.

\begin{theo}\label{filtration}
	For every $G \in \mathfrak{f}$,
\[ {\mathbb{K}}^{- \infty}(\mathcal{B} (CX,CY \cup X,p_X;R))^G \simeq {\mathbb{K}}^{- \infty}(\mathcal{B} (\Sigma X,\Sigma Y,p_X;R))^G. \]
\end{theo}

	Given $G \in \mathfrak{f}$, consider the subgroup lattice for $G$, and let $H\leq G$. Define the distance from $G$ to $H$, $dist(H)$, to be the maximum number of steps needed to reach $H$ from $G$ on the subgroup lattice. Notice that $dist(gHg^{-1})=dist(H)$ for each $g\in G$.

	Let $n=dist(1)$, and let $l_i$ be the number of conjugacy classes of subgroups with distance $i$ from $G$. For each $i$, $1\leq i\leq n-1$, choose a representative, $H_{i,j}$, $1\leq j\leq l_i$, for each of the conjugacy classes with distance $i$ from $G$. Order the $H_{i,j}$'s using the dictionary order on the indexing set. Now re-index the sequence of subgroups according to the ordering so that
\[H_1=H_{1,1} \ \ , \ \ H_2=H_{1,2} \ \ , \ ... \ , \ \ H_m=H_{n-1,l_{n-1}}. \]
For convenience define $H_0=G$ and $H_{m+1}=1$.

	For each $j$, $0\leq j\leq m+1$, define $C_j=\bigcup_{g \in G} X^{gH_jg^{-1}}$, which will often be referred to as a {\it conjugacy class of fixed sets}.

	For each $k$, $0\leq k\leq m+1$, define $Z_k=\bigcup_{0 \leq j \leq k} C_j$. Also define $Y_k=Z_k \cap Y$. Since it is possible that $Y_k= \emptyset$, we define $C(\emptyset)=\{ 0 \}$ and $\Sigma (\emptyset)=\{ 0,1 \}$ for the convenience of notation.
	
	Notice that $Z_k$ is contractible for every $k$, $0\leq k\leq m+1$. Also notice that since $gX^H=X^{gHg^{-1}}$, each of the conjugacy classes of fixed sets is $G$-invariant. Therefore $Z_k$ is $G$-invariant for every $k$, $0\leq k\leq m+1$.

\begin{rema}\label{remark1}
	Let $z\in Z_k-Z_{k-1}$, $1\leq k \leq m$. Then the stabilizer of $z$ is a conjugate of $H_k$.
\end{rema}

\begin{proof}
	Certainly $z\in X^{gH_kg^{-1}}$ for some $g\in G$. The only question is whether there is a larger subgroup fixing $z$. Suppose that the stabilizer of $z$ is $H$. Then $gH_kg^{-1}$ is a proper subgroup of $H$. Hence, $dist(H)<dist(gH_kg^{-1})$. Therefore $X^H$ appears in a previous conjugacy class of fixed sets. But this implies $z\in Z_{k-1}$.
\end{proof}

\begin{rema}\label{remark2}
	Let $H,K \leq G$, $H\neq K$ with $X^H,X^K \subseteq Z_k$ for some $k$, $1\leq k \leq m$. Then $X^H \cap X^K \subseteq Z_{k-1}$.
\end{rema}

\begin{proof}
	Let $\langle H,K \rangle$ denote the smallest subgroup of $G$ containing both $H$ and $K$. Then $X^H \cap X^K=X^{\langle H,K \rangle}$. Since $H\neq K$, at least one of the them must be a proper subgroup of $\langle H,K \rangle$. Therefore $X^{\langle H,K \rangle}$ appears in one of the previous conjugacy classes of fixed sets.
\end{proof}

	To simplify the notation set
\begin{align}
	&{\mathcal{B}}_G (C(Z_k);R)= \mathcal{B} (C(Z_k),C(Y_k) \cup Z_k,p_{Z_k};R)^G; \notag\\ 
	&{\mathcal{B}}_G (\Sigma (Z_k);R) = \mathcal{B} (\Sigma (Z_k),\Sigma (Y_k),p_{Z_k};R)^G; \notag\\
	&{\mathcal{B}}_G (C(Z_k);R)_{Z_{k-1}} = {\mathcal{B}}_G (C(Z_k);R)_{C(Y_{k-1}) \cup Z_{k-1}}; \notag\\
        &{\mathcal{B}}_G (\Sigma (Z_k);R)_{Z_{k-1}} = {\mathcal{B}}_G (\Sigma(Z_k);R)_{\Sigma (Y_{k-1})}; \notag\\
	&{\mathcal{B}}_G (C(Z_k);R)^{>Z_{k-1}} = {\mathcal{B}}_G (C(Z_k);R)^{(C(Y_k) \cup Z_k)-(C(Y_{k-1}) \cup Z_{k-1})}; \notag\\
	&{\mathcal{B}}_G (\Sigma (Z_k);R)^{>Z_{k-1}} = {\mathcal{B}}_G (\Sigma(Z_k);R)^{\Sigma(Y_k)-\Sigma (Y_{k-1})} . \notag
\end{align}

	Theorem $\ref{filtration}$ is proved by induction on the chain \[X^G=Z_0 \subseteq Z_1 \subseteq \dddot{} \subseteq Z_m \subseteq Z_{m+1}=X. \]
This is accomplished with the next two lemmas.

\begin{lem}\label{lemma2}
	For each $k$, $1\leq k\leq m+1$,
\begin{enumerate}
\item[1.] ${\mathcal{B}}_G (C(Z_k);R)^{>Z_{k-1}} \cong \mathcal{B} (C({Z_k} / {G});R[H_k])^{>Z_{k-1} / G}$.
\item[2.] ${\mathcal{B}}_G (\Sigma(Z_k);R)^{>Z_{k-1}} \cong \mathcal{B} (\Sigma({Z_k} / {G});R[H_k])^{>Z_{k-1} / G}$.
\end{enumerate}
\end{lem}

\begin{proof}
	Since we are taking germs away from $C(Y_{k-1}) \cup Z_{k-1}$ and $C(Y_{k-1} / G) \cup Z_{k-1} / G$, every morphism has a representative that is zero on $C(Z_{k-1})$ and on $C(Z_{k-1} / G)$, respectively. It is therefore irrelevant what the objects over $C(Z_{k-1})$ and $C(Z_{k-1} / G)$ are. Nevertheless, we want to establish a one-one correspondence between the objects of the two categories. Given an object $A$ in ${\mathcal{B}}_G (C(Z_k);R)^{>Z_{k-1}}$, we define an object $A'$ in $\mathcal{B} (C(Z_k / G);R[H_k])^{>Z_{k-1} / G}$ as follows: $A'_{[x]}:=A_x$ if $[x]\notin C(Z_{k-1} / G)$; otherwise $A'_{[x]}$ is defined to be a free $R[H_k]$-module of rank equal to ${\rm rank}_{RG}\big(\bigoplus_{x'\in [x]} A_{x'}\big)$. If $B$ is an object in $\mathcal{B} (C({Z_k} / {G});R[H_k])^{>Z_{k-1} / G}$, then an object $\bar{B}$ in ${\mathcal{B}}_G (C(Z_k);R)^{>Z_{k-1}}$ is defined as follows: $\bar{B}_x:=B_{[x]}$ if $x\notin C(Z_{k-1})$; otherwise $\bigoplus_{x'\in [x]} \bar{B}_{x'}$ is defined to be a free $RG$-module of rank equal to ${\rm rank}_{R[H_k]}(B_{[x]})$. This gives the desired one-one correspondence.

	Since we are taking germs, the components of a morphism need to become small. Therefore, using Remarks $\ref{remark1}$ and $\ref{remark2}$, and the fact that we are working with a finite group action, we can be sure that non-zero components of a morphism have the same isotropy, namely a conjugate of $H_k$. Therefore by equivariance, there is only one choice when lifting a morphism. This completes the proof of the first part of this lemma. The same proof works for the second part, replacing cones with suspensions.
\end{proof}

\begin{lem}\label{lemma3}
	For each $k$, $1\leq k\leq m+1$, 
\[ {\mathbb{K}}^{- \infty}\left(\mathcal{B} (C({Z_k} / G);R[H_k])^{>Z_{k-1} / G}\right) \simeq {\mathbb{K}}^{- \infty}\left(\mathcal{B} (\Sigma({Z_k} / {G});R[H_k])^{>Z_{k-1} / G}\right). \]
\end{lem}

\begin{proof}
	Let $A=Z_{k-1} / G$. Consider the following commutative diagram:
\[ \xymatrix{
		\mathcal{B}(C({Z_k} / G);R[H_k])_A \ar[d]_a \ar[r] & \mathcal{B}(C({Z_k} / G);R[H_k]) \ar[d]_b \ar[r] & \mathcal{B}(C({Z_k} / G);R[H_k])^{>A} \ar[d]_c \\
		\mathcal{B}(\Sigma({Z_k} / {G});R[H_k])_A \ar[r] & \mathcal{B}(\Sigma({Z_k} / {G});R[H_k]) \ar[r] & \mathcal{B}(\Sigma({Z_k} / {G});R[H_k])^{>A}. } \]

	Since $Z_k$ is contractible, the reduced Steenrod homology of ${Z_k} / G$ is trivial by Corollary~\ref{zero}. Therefore $b$ induces a weak homotopy equivalence by Proposition $\ref{key}$.
	
	By Lemma $\ref{p_x}$,
\[ \mathcal{B} (C(Z_k / G);R[H_k])_{Z_{k-1} / G} \cong \mathcal{B} (C(Z_{k-1} / G);R[H_k]) \]
and
\[ \mathcal{B} (\Sigma(Z_k / G);R[H_k])_{Z_{k-1} / G} \cong \mathcal{B} (\Sigma(Z_{k-1} / G);R[H_k]). \]
Since $Z_{k-1}$ is contractible, the reduced Steenrod homology of ${Z_{k-1}} / G$ is trivial by Corollary~\ref{zero}. Thus, by Proposition $\ref{key}$, $a$ induces a weak homotopy equivalence. Since each row in the diagram is a Karoubi filtration, the Five Lemma shows that $c$ also induces a weak homotopy equivalence.
\end{proof}

We are now ready to prove Theorem $\ref{filtration}$.

\begin{proof}[Proof of Theorem \ref{filtration}]
	From the definition of a fixed category, ${\mathcal{B}}_G (C(X^G);R) \cong \mathcal{B} (C(X^G);RG)$ and $\mathcal{B} (\Sigma(X^G);RG) \cong {\mathcal{B}}_G (\Sigma(X^G);R)$. Since $X^G$ is contractible, Proposition~\ref{key} implies ${\mathbb{K}}^{- \infty}(\mathcal{B} (C(X^G);RG)) \simeq {\mathbb{K}}^{- \infty}(\mathcal{B} (\Sigma(X^G);RG))$. Therefore
\[ {\mathbb{K}}^{- \infty}\left({\mathcal{B}}_G (C(X^G);R)\right) \simeq {\mathbb{K}}^{- \infty}\left({\mathcal{B}}_G (\Sigma(X^G);R)\right). \]
This completes the base case of the induction. 
	
	Assume now that ${\mathbb{K}}^{- \infty}({\mathcal{B}}_G (C(Z_{k-1});R)) \simeq {\mathbb{K}}^{- \infty}({\mathcal{B}}_G (\Sigma(Z_{k-1});R))$. We want to show that ${\mathbb{K}}^{- \infty}({\mathcal{B}}_G (C(Z_k);R)) \simeq {\mathbb{K}}^{- \infty}({\mathcal{B}}_G (\Sigma(Z_k);R))$.
Consider the following commutative diagram:
\[ \xymatrix{
	{\mathcal{B}}_G (C(Z_k);R)_{Z_{k-1}} \ar[d]_a \ar[r] & {\mathcal{B}}_G (C(Z_k);R) \ar[d]_b \ar[r] & {\mathcal{B}}_G (C(Z_k);R)^{>Z_{k-1}} \ar[d]_c \\
	{\mathcal{B}}_G (\Sigma(Z_k);R)_{Z_{k-1}} \ar[r] & {\mathcal{B}}_G (\Sigma(Z_k);R) \ar[r] & {\mathcal{B}}_G (\Sigma(Z_k);R)^{>Z_{k-1}}. } \]

	By Lemma $\ref{p_x}$, Proposition $\ref{key}$, and the induction hypothesis, $a$ induces a weak homotopy equivalence. By Lemmas $\ref{lemma2}$ and $\ref{lemma3}$, $c$ induces a weak homotopy equivalence. Since each row in the diagram is a Karoubi filtration, the Five Lemma shows that $b$ also induces a weak homotopy equivalence.
\end{proof}

This completes the proof of Theorem~\ref{main}.

\section{A Davis-L\"{u}ck Spectral Sequence}\label{dlss}

	In this section we develop a spectral sequence that comes out of the Davis-L\"{u}ck machinery introduced in~\cite{dl}. We begin by recalling some terminology. The {\it orbit category} Or$(\Gamma)$ has as objects homogeneous $\Gamma$-spaces, $\Gamma / H$, and as morphisms $\Gamma$-maps. An Or$(\Gamma)$-space is a functor from Or$(\Gamma)$ to the category of spaces. Similarly, an Or$(\Gamma)$-spectrum is a functor from Or$(\Gamma)$ to the category of spectra. A map of Or$(\Gamma)$-spaces is a natural transformation. A map of Or$(\Gamma)$-spaces is called a {\it weak homotopy equivalence} if it is a weak homotopy equivalence evaluated at every $\Gamma / H$. If $\mathscr{X}$ is a contravariant Or$(\Gamma)$-space and $\mathscr{Y}$ is a covariant Or$(\Gamma)$-space, then their tensor product is the space
\[ \mathscr{X} \otimes_{{\rm Or}(\Gamma)} \mathscr{Y} = \coprod_{\Gamma / H} \mathscr{X}(\Gamma / H) \times \mathscr{Y}(\Gamma / H) / \sim, \]
where $\sim$ is the equivalence relation generated by $(\mathscr{X}(\phi)(x),y) \sim (x,\mathscr{Y}(\phi)(y))$ for all morphisms $\phi: \Gamma / H \to \Gamma /K$ in Or$(\Gamma)$, and points $x \in \mathscr{X}(\Gamma / K)$ and $y \in \mathscr{Y}(\Gamma / H)$. This notion extends to pointed Or$(\Gamma)$-spaces by replacing disjoint unions and Cartesian products with wedge products and smash products. This construction can then be generalized to obtain the tensor product spectrum of a contravariant Or$(\Gamma)$-space with a covariant Or$(\Gamma)$-spectrum. Homology theories can now be defined. If $\mathbb{E}$ is a covariant  Or$(\Gamma)$-spectrum and $(\mathscr{X},\mathscr{A})$ is a pair of  Or$(\Gamma)$-spaces, then
\[ H^{{\rm Or}(\Gamma)}_{n}(\mathscr{X},\mathscr{A};\mathbb{E})=\pi_n(\mathscr{X}_+ \cup_{\mathscr{A}_+} C(\mathscr{A}_+) \otimes_{{\rm Or}(\Gamma)} \mathbb{E}), \]
where $C(\mathscr{A}_+)$ denotes the reduced cone of $\mathscr{A}_+$.
	
	Let $X$ be a $\Gamma$-CW complex with $n$-skeleton $X_n$ and let $\Sigma_n$ be the set of $n$-cells in $X / \Gamma$. Assume that if a cell in $X$ is fixed by an element of $\Gamma$ then it is pointwise fixed. Let $\mathcal{F}$ be the family of subgroups of $\Gamma$ that appear as stabilizers of cells in $X$. For each $\sigma \in \Sigma_n$ choose an $n$-cell in $X$, $\tilde{\sigma}$, representing $\sigma$.

	Let $\mathscr{X}=\Map_{\Gamma}(-,X)$ and $\mathscr{X}_n=\Map_{\Gamma}(-,X_n)$. Then $\mathscr{X}$ and $\mathscr{X}_n$ are contravariant free Or$(\Gamma)$-CW complexes. Let $\mathbb{E}$ be a covariant Or$(\Gamma)$-spectrum. There is a spectral sequence converging to $H^{{\rm Or}(\Gamma)}_{p+q}(\mathscr{X};\mathbb{E})$ whose $E^1$-term is given by $E^1_{p,q}=H^{{\rm Or}(\Gamma)}_{p+q}(\mathscr{X}_p,\mathscr{X}_{p-1};\mathbb{E})$ and whose first differential, $d^1_{p,q}$, is the composition 
\[ H^{{\rm Or}(\Gamma)}_{p+q}(\mathscr{X}_p,\mathscr{X}_{p-1};\mathbb{E}) \to  H^{{\rm Or}(\Gamma)}_{p+q-1}(\mathscr{X}_{p-1};\mathbb{E}) \to   H^{{\rm Or}(\Gamma)}_{p+q-1}(\mathscr{X}_{p-1},\mathscr{X}_{p-2};\mathbb{E}), \]
where the first map is the boundary operator of the pair $(\mathscr{X}_p,\mathscr{X}_{p-1})$ and the second is induced by inclusion~\cite[Theorem 4.7]{dl}. We wish to analyze this spectral sequence more closely.

	Since we are assuming that cells are pointwise fixed, $X^H$ is a CW complex. Therefore the natural transformation $\mathscr{X}_{p_+} \negthinspace \cup_{\mathscr{X}_{{p-1}_+}} \negmedspace C(\mathscr{X}_{{p-1}_+}) \to \mathscr{X}_{p_+} / \mathscr{X}_{{p-1}_+}$ is a weak homotopy equivalence of free Or$(\Gamma)$-CW complexes and hence a homotopy equivalence of Or$(\Gamma)$-spaces~\cite[Corollary 3.5]{dl}. Therefore $H^{{\rm Or}(\Gamma)}_{p+q}(\mathscr{X}_p,\mathscr{X}_{p-1};\mathbb{E}) = \pi_{p+q}(\mathscr{X}_{p_+} \negthinspace \cup_{\mathscr{X}_{{p-1}_+}} \negmedspace C(\mathscr{X}_{{p-1}_+}) \otimes_{{\rm Or}(\Gamma)} \mathbb{E} ) \cong \pi_{p+q}(\mathscr{X}_{p_+} / \mathscr{X}_{{p-1}_+} \negthinspace \otimes_{{\rm Or}(\Gamma)} \mathbb{E} )$. Consider 
\[ \mathscr{X}_{p_+} / \mathscr{X}_{{p-1}_+} \negthinspace \otimes_{{\rm Or}(\Gamma)} \mathbb{E}_n = {\bigvee_{H \in \mathcal{F}} {X^H_{p_+}} / {X^H_{{p-1}_+}} \negthinspace \wedge \mathbb{E}_n(\Gamma / H)} / \sim. \]
Let $x$ be in the interior of a $p$-cell, $\tau \subset X^H_p$, and $s \in \mathbb{E}_n(\Gamma / H)$. Then there is a $p$-cell, $\tilde{\sigma}$, with $\tau=\gamma\tilde{\sigma}$ for some $\gamma\in \Gamma$. This implies that $H \leq \Gamma_{\tau}$ and $\Gamma_{\tau}=\Gamma_{\gamma\tilde{\sigma}}=\gamma\Gamma_{\tilde{\sigma}}\gamma^{-1}$. So $\phi:\Gamma / H \to \Gamma / \Gamma_{\tilde{\sigma}}$, defined by $\phi(H)=\gamma\Gamma_{\tilde{\sigma}}$, is a morphism in Or$(\Gamma)$. Therefore $h(x,s)=(\phi^{*^{-1}}(x),\phi_*(s))$ well-defines a map of spectra
\[ \mathscr{X}_{p_+} / \mathscr{X}_{{p-1}_+} \negthinspace \otimes_{{\rm Or}(\Gamma)} \mathbb{E} \to \bigvee_{\sigma \in \Sigma_p} \tilde{\sigma}_+ / \partial\tilde{\sigma}_+ \negthinspace \wedge \mathbb{E}(\Gamma / \Gamma_{\tilde{\sigma}}), \]
that is an inverse to the inclusion map, and hence an equivalence of spectra. This implies that
\[ E^1_{p,q} \cong \bigoplus_{\sigma \in {\Sigma}_p} \pi_{p+q}\big(S^p \wedge \mathbb{E}(\Gamma / \Gamma_{\tilde{\sigma}})\big)  \cong \bigoplus_{\sigma \in {\Sigma}_p} \pi_q\big(\mathbb{E}(\Gamma / \Gamma_{\tilde{\sigma}})\big). \]
	
	As one would expect, the first differential is calculated by generalizing the derivation of the boundary operator in cellular homology. Specifically, we generalize the proof of~\cite[Proposition 10.11]{switzer}. 
	
	Let $\alpha \in \Sigma_p$ and $\beta \in \Sigma_{p-1}$. Let $X_{p-1}^{\beta}$ be the set of all points in $X$ that do not lie in the interior of the $\Gamma$-cell represented by $\tilde{\beta}$. Let $\bar{f}_{\tilde{\alpha}}:S^{p-1} \to X_{p-1} / X_{p-1}^{\beta}$ and $\bar{g}:\bigvee_{\lambda \in \Gamma / \Gamma_{\tilde{\beta}}} S^{p-1}_{\lambda} \to X_{p-1} / X_{p-1}^{\beta}$ be the maps induced by the attaching map for $\tilde{\alpha}$ and the characteristic map for the $\Gamma$-cell represented by $\tilde{\beta}$, respectively. Let $h=\bar{g}^{-1} \circ \bar{f}_{\tilde{\alpha}}$ and $h_{\tilde{\alpha},\gamma \tilde{\beta}}=p_{\gamma \tilde{\beta}} \circ h$, where $p_{\gamma \tilde{\beta}}$ is projection onto the $\gamma \Gamma_{\tilde{\beta}}$ factor of $\bigvee_{\lambda \in \Gamma / \Gamma_{\tilde{\beta}}} S^{p-1}_{\lambda}$.
	
	Using a diagram similar to the one used in~\cite[Proposition 10.11]{switzer}, we see that the restriction of $d^1_{p,q}$ to $\pi_q\big(\mathbb{E}(\Gamma / \Gamma_{\tilde{\alpha}})\big) \to \pi_q\big(\mathbb{E}(\Gamma / \Gamma_{\tilde{\beta}})\big)$ is induced by a map
\[ d: S^{p-1} \wedge \mathbb{E}(\Gamma / \Gamma_{\tilde{\alpha}}) \rightarrow S^{p-1} \wedge \mathbb{E}(\Gamma / \Gamma_{\tilde{\beta}}), \]
which we study in the diagram below. Let $I_{\beta}$ be an indexing set for the number of translates of $\tilde{\beta}$ that are contained in the boundary of $\tilde{\alpha}$. This set must be finite. If $I_{\beta}$ is empty then $d$ maps everything to the basepoint and therefore induces the zero map on homotopy groups. Assume that $I_{\beta}$ is non-empty. Choose representatives $\gamma_i \in \Gamma$, one for each $i \in I_{\beta}$, such that $\gamma_i \tilde{\beta} \subset \partial \tilde{\alpha}$. The above map then fits into the following commutative diagram

\[ \xymatrix{
	\Big(\bigvee_i S^{p-1}_i \Big) \wedge \mathbb{E}(\Gamma / \Gamma_{\tilde{\alpha}}) \ar[r]^{\cong} & \bigvee_i S^{p-1}_i \wedge \mathbb{E}(\Gamma / \Gamma_{\tilde{\alpha}}) \ar[r]^{\bigvee 1 \wedge \psi_{i_*}} &	\bigvee_i S^{p-1}_i \wedge \mathbb{E}(\Gamma / \Gamma_{\tilde{\beta}}) \ar[dl]^{{\rm proj.}}\\
	S^{p-1} \wedge \mathbb{E}(\Gamma / \Gamma_{\tilde{\alpha}}) \ar[u]^{h \wedge 1} \ar[r]^{d} & S^{p-1} \wedge \mathbb{E}(\Gamma / \Gamma_{\tilde{\beta}}), } \]
where $\psi_i:\Gamma / \Gamma_{\tilde{\alpha}} \to \Gamma / \Gamma_{\tilde{\beta}}$ is the Or$(\Gamma)$-morphism defined by $\psi_i(\Gamma_{\tilde{\alpha}})=\gamma_i\Gamma_{\tilde{\beta}}$. Note that $\psi_i$ is independent of the choice of representative. This implies that $\sum_{i \in I_{\beta}}({h_{\tilde{\alpha},\gamma_i \tilde{\beta}}}_* \wedge {\psi_i}_*)=d_*: \pi_{p+q-1}\big(S^{p-1} \wedge \mathbb{E}(\Gamma / \Gamma_{\tilde{\alpha}})\big) \to \pi_{p+q-1}\big(S^{p-1} \wedge \mathbb{E}(\Gamma / \Gamma_{\tilde{\beta}})\big)$. Since ${h_{\tilde{\alpha},\gamma_i \tilde{\beta}}}_*$ is just multiplication by ${\rm deg}(h_{\tilde{\alpha},\gamma_i \tilde{\beta}})=[\tilde{\alpha}:\gamma_i\tilde{\beta}]$, we see that the restriction of $d^1_{p,q}$ to $\pi_q\big(\mathbb{E}(\Gamma / \Gamma_{\tilde{\alpha}})\big) \to \pi_q\big(\mathbb{E}(\Gamma / \Gamma_{\tilde{\beta}})\big)$ is $\sum_{i \in I_{\beta}} [\tilde{\alpha}:\gamma_i\tilde{\beta}]\psi_{i_*}$. We have proven the following.

	\begin{theo}\label{spectral}
	There is a spectral sequence converging to $H^{{\rm Or}(\Gamma)}_{p+q}(\mathscr{X};\mathbb{E})$ whose $E^1$-term is given by 
\[ E^1_{p,q}=\bigoplus_{\sigma \in {\Sigma}_p} \pi_q\big(\mathbb{E}(\Gamma / \Gamma_{\tilde{\sigma}})\big) \]
and whose first differential, $d^1_{p,q}:E^1_{p,q} \to E^1_{p-1,q}$, is given by
\[ d^1_{p,q}(a)=\sum_{\beta \in \Sigma_{p-1}}\Big(\sum_{i \in I_{\beta}} [\tilde{\alpha}:\gamma_i\tilde{\beta}]\psi_{i_*}(a)\Big), \]
where $a \in \pi_q\big(\mathbb{E}(\Gamma / \Gamma_{\tilde{\alpha}})\big)$.
	\end{theo}

	This spectral sequence is a generalization of the one given in \cite[Chapter VII]{brown}. In the next section we will use it with the algebraic $K$-theory Or$(\Gamma)$-spectrum, $\mathbb{K}^{alg}_R$, associated to a chosen ring $R$. Recall that $\pi_q\big(\mathbb{K}^{alg}_R(\Gamma / H)\big) \cong K_q(RH)$.

\section{Virtually Polycyclic Groups}\label{vpg}

A group $\Gamma$ is called {\it virtually polycyclic} if there is a finite sequence of normal subgroups
\[ 1=N_1 \triangleleft \dotsb \triangleleft N_k=\Gamma \]
such that $N_{i+1} / N_i$ is either infinite cyclic or finite. In~\cite[Theorem 1]{moody}, Moody proved that if $\Gamma$ is a virtually polycyclic group, then the induction map
\[ \bigoplus_{H \in \mathfrak{f}} G_0(RH) \to G_0(R\Gamma) \]
is a surjection. If $R$ is a regular ring then $G_0(R) \cong K_0(R)$. This is the case for $k\Gamma$, when $k$ is a field of characteristic zero and $\Gamma$ is a virtually polycyclic group~\cite{bass}. We have the following proposition as a special case of Moody's result.

	\begin{prop}\label{surjection}
	Let $\Gamma$ be a virtually polycyclic group and $k$ be a field of characteristic zero. Then the induction map
\[ \colim_{H \in \mathfrak{f}} K_0(kH) \to K_0(k\Gamma) \]
is a surjection.
	\end{prop}

	In this section we will use Theorems~\ref{main} and $\ref{spectral}$ to show that this map is in fact an isomorphism. 	
		
	\begin{prop}\label{h_0}
	Let $\Gamma$ be a group and $k$ be a field of characteristic zero. Let $\mathscr{E}=\Map_{\Gamma}(-,E\Gamma(\mathfrak{f}))$. Then
\[ H_0^{{\rm Or}(\Gamma)} (\mathscr{E};\mathbb{K}^{alg}_k) \cong \colim_{H \in \mathfrak{f}} K_0(kH). \]
	\end{prop}
	
	\begin{proof}
	  Choose $E\Gamma(\mathfrak{f})$ so that if a cell is fixed by a group element, then it is pointwise fixed by that element~\cite[Appendix]{fj}. Consider the spectral sequence in Theorem~\ref{spectral} with $\mathbb{E}=\mathbb{K}^{alg}_k$ and $\mathscr{X}=\mathscr{E}$. Since $k$ is a field of characteristic zero and $H$ is a finite group, the negative $K$-groups of $kH$ are all zero~\cite{bass}. Therefore the spectral sequence collapses at $E^2$. In particular, $H_0^{{\rm Or}(\Gamma)} (\mathscr{E};\mathbb{K}^{alg}_k) \cong E^2_{0,0}=E^1_{0,0} / {\rm im}(d^1_{1,0})$. We want to show that $E^2_{0,0} \cong \colim_{H \in \mathfrak{f}} K_0(kH)$. The map
\[ E^1_{0,0} \cong \bigoplus_{x \in {\Sigma}_0} K_q(k\Gamma_{\tilde{x}}) \hookrightarrow \bigoplus_{H \in \mathfrak{f}} K_q(kH) \xrightarrow{\text {\it q}} \colim_{H \in \mathfrak{f}} K_0(kH) \]
induces a map $E^2_{0,0} \to \colim_{H \in \mathfrak{f}} K_0(kH)$ since the image of 
\[ d^1_{1,0}: \bigoplus_{\sigma \in {\Sigma}_1} K_q(k\Gamma_{\tilde{\sigma}}) \to \bigoplus_{x \in {\Sigma}_0} K_q(k\Gamma_{\tilde{x}}) \]
is contained in the kernel of $q$. 
Since inner automorphisms induce the identity on $K$-theory, we need only consider one subgroup from each conjugacy class in the colimit. Furthermore, only the maximal finite subgroups of $\Gamma$ are needed in the colimit. Therefore this map is a surjection. To show that this map is an isomorphism we need to show that the set of relations that generate the kernel of $q$ are contained in the image of $d^1_{1,0}$.

	It suffices to prove that if $G_1$ and $G_2$ are maximal finite subgroups of $\Gamma$, and $a_1 \in K_0(kG_1)$ is identified with $a_2 \in K_0(kG_2)$ in the colimit, then $a_1-a_2$ is in the image of $d^1_{1,0}$. Since $G_1$ and $G_2$ are maximal, they must appear as stabilizers of $0$-cells in $E\Gamma(\mathfrak{f})$. Since we only need one group for each conjugacy class, we assume that $G_1=\Gamma_{\tilde{x}}$ and $G_2=\Gamma_{\tilde{y}}$ for the chosen representatives $\tilde{x}$, $\tilde{y} \in E\Gamma(\mathfrak{f})$, for some $x$, $y \in \Sigma_0$. Since $a_1$ is identified with $a_2$, there exists an $H\in \mathfrak{f}$ with $H \leq G_1$ and $H \leq G_2$, and an $a\in K_0(kH)$, such that $a_1$ and $a_2$ are the image of $a$ under the homomorphisms $K_0(kH) \to K_0(kG_1)$ and $K_0(kH) \to K_0(kG_2)$, respectively. Since $E\Gamma(\mathfrak{f})^H$ is contractible, there is a path of $1$-cells in $E\Gamma(\mathfrak{f})^H$ connecting $\tilde{x}$ and $\tilde{y}$. Let $\tau_1,...,\tau_n$ be a sequence of $1$-cells in $E\Gamma(\mathfrak{f})^H$ connecting $\tilde{x}$ to $\tilde{y}$, with $\partial_0 \tau_1=\tilde{x}$, $\partial_1 \tau_n=\tilde{y}$, and $\partial_0 \tau_{i+1}=\tilde{x}_i=\partial_1 \tau_i$ for $1\leq i \leq n-1$. Note that $H \leq \Gamma_{\tau_i}$ and $H \leq \Gamma_{x_i}$ for every $i$. Let $\tilde{\sigma}_i$ be the chosen representative for the orbit containing $\tau_i$ and let $\tilde{z}_i$ be the chosen representative for the orbit containing $x_i$. For the convenience of notation, let $\tilde{z}_0=\tilde{x}$ and $\tilde{z}_n=\tilde{y}$. Note that $\Gamma_{\tilde{\sigma}_i}$ is conjugate to $\Gamma_{\tau_i}$, and that $\Gamma_{\tilde{z}_i}$ is conjugate to $\Gamma_{x_i}$. Let $a_{\tilde{\sigma}_i}$ be the image of $a$ under the composition $K_0(kH) \to K_0(k\Gamma_{\tau_i}) \to K_0(k\Gamma_{\tilde{\sigma}_i})$, and let $a_{\tilde{z}_i}$ be the image of $a$ under the composition $K_0(kH) \to K_0(k\Gamma_{x_i}) \to K_0(k\Gamma_{\tilde{z}_i})$. Then $d^1_{1,0}(a_{\tilde{\sigma}_i})=a_{\tilde{z}_i}-a_{\tilde{z}_{i-1}}$ for $1\leq i \leq n$. Therefore
\[d^1_{1,0}\Big(\sum_{i=1}^n a_{\tilde{\sigma}_i} \Big)=\sum_{i=1}^n (a_{\tilde{z}_i}-a_{\tilde{z}_{i-1}})=a_{\tilde{x}}-a_{\tilde{y}}=a_1-a_2.  \]

Finally, we remark that although there may be more than one choice for $E\Gamma(\mathfrak{f})$, the result does not depend on this choice since any two models are $\Gamma$-equivariantly homotopy equivalent.
	\end{proof}

	By comparing the $E^1_{0,0}$-term for the left and right hand sides of the Davis-L\"{u}ck assembly map, we see that the assembly map in dimension zero,
\[ \colim_{H \in \mathfrak{f}} K_0(kH) \to K_0(k\Gamma), \]
is just the induction map. Since the Davis-L\"{u}ck assembly map is homotopy equivalent to the continuously controlled assembly map \cite[Theorem 8.3]{hp}, we can apply Theorem~\ref{main} to show that this map is injective. Before we do this we need to know that there is a model for $E\Gamma(\mathfrak{f})$ that satisfies all of the necessary conditions.

 In~\cite[Theorem 3]{wilking}, Wilking proved that a group $\Gamma'$ is isomorphic to a discrete, cocompact subgroup of a semidirect product $S_1 \rtimes K$, where $S_1$ is a connected, simply connected solvable Lie group and $K$ is a compact subgroup of the automorphism group ${\rm Aut}(S_1)$ if and only if $\Gamma'$ is virtually polycyclic containing no non-trivial finite normal subgroups. This implies that $S_1 \rtimes K / K \cong \mathbb{R}^n$ is an $E\Gamma'(\mathfrak{f})$. By compactifying with an $(n-1)$-sphere at infinity we see that Theorem~\ref{main} applies to such $\Gamma'$.
	
	\begin{lem}\label{nicas}
	  Suppose that $0 \to N \to \Gamma \xrightarrow{\text {\it p}} \Gamma' \to 0$ is a short exact sequence, where $N$ is a finite normal subgroup of $\Gamma$. If $E$ is a model for $E\Gamma'(\mathfrak{f})$, then $E$ is also a model for $E\Gamma(\mathfrak{f})$, where $\Gamma$ acts via the homomorphism $p$.
	\end{lem}
	
	\begin{proof}
	  Let $H$ be a subgroup of $\Gamma$. We need to show that $E^H=E^{p(H)}$ is contractible if $H$ is finite, and that it is empty otherwise. This is true since $p(H)$ is finite if and only if $H$ is finite.
	\end{proof}
	
	All virtually polycyclic groups contain a torsion-free subgroup of finite index~\cite{mcconnell}. Therefore if $H$ is a finite normal subgroup of a virtually polycyclic group $\Gamma$, then it must be contained in a maximal finite normal subgroup, $N$. This implies that $\Gamma'=\Gamma / N$ is a virtually polycyclic group with no non-trivial finite normal subgroups. Now using Wilking's result, mentioned above, and Lemma~\ref{nicas}, the conditions of Theorem~\ref{main} are satisfied by all virtually polycyclic groups. Therefore we have proven:

	\begin{coro}\label{iso}
	Let $\Gamma$ be a virtually polycyclic group and $k$ a field of characteristic zero. Then the induction map
\[ \colim_{H \in \mathfrak{f}} K_0(kH) \to K_0(k\Gamma) \]
is an isomorphism.
	\end{coro}

\end{document}